\documentclass[11pt]{article}
\usepackage{amsmath, amsthm, amssymb}    
\usepackage{bridges}			
\usepackage{graphicx}			
\usepackage[colorlinks=true, urlcolor=blue, citecolor=black, linkcolor=black]{hyperref}  
\usepackage{subcaption}			
\usepackage{nicefrac}

\newcommand{\R}{\mathbb{R}}
\newcommand{\C}{\mathbb{C}}
\newcommand{\I}{\mathbf{i}}
\newcommand{\J}{\mathbf{j}}
\newcommand{\K}{\mathbf{k}}

\newcommand{\St}{\operatorname{St}}
\newcommand{\Gr}{\operatorname{Gr}}

\title{Stiefel Manifolds and Polygons}
\author{Clayton Shonkwiler\\
\vspace{10pt}\\
Department of Mathematics, Colorado State University; clay@shonkwiler.org}

\date{}				

\begin{document}

\maketitle

\thispagestyle{empty}

\begin{abstract}
	Polygons are compound geometric objects, but when trying to understand the expected behavior of a large collection of random polygons -- or even to formalize what a random polygon is -- it is convenient to interpret each polygon as a point in some parameter space, essentially trading the complexity of the object for the complexity of the space. In this paper I describe such an interpretation where the parameter space is an abstract but very nice space called a \emph{Stiefel manifold} and show how to exploit the geometry of the Stiefel manifold both to generate random polygons and to morph one polygon into another.
\end{abstract}



\section*{Introduction} 
The goal is to describe an identification of polygons with points in certain nice parameter spaces. For example, every triangle in the plane can be identified with a pair of orthonormal vectors in 3-dimensional space~$\R^3$, and hence planar triangle shapes correspond to points in the space $\St_2(\R^3)$ of all such pairs, which is an example of a Stiefel manifold. While this space is defined somewhat abstractly and is hard to visualize -- for example, it cannot be embedded in Euclidean space of fewer than 5 dimensions~\cite{Mahowald:1962ep} -- it is easy to visualize and manipulate points in it: they're just pairs of perpendicular unit vectors. Moreover, this is a familiar and well-understood space in the setting of differential geometry. For example, the group $SO(3)$ of rotations of 3-space acts transitively on it, meaning that there is a rotation of space which transforms any triangle into any other triangle.

This identification between triangles and pairs of orthonormal vectors generalizes: we will represent an $n$-gon in the plane by a pair of orthonormal vectors in $\R^n$. The space $\St_2(\R^n)$ of such pairs is again a Stiefel manifold -- in general, the \emph{Stiefel manifold} $\St_k(\R^n)$ is the collection of all $k$-tuples of mutually orthonormal vectors in $\R^n$ -- and again there is a natural transitive action, now by the group $SO(n)$ of all rotations of~$\R^n$. The existence of transitive group actions is special in both probability theory and differential geometry, and in particular means that the Stiefel manifold (and, hence, the collection of planar $n$-gons) has a natural probability measure (\emph{Haar measure}) which is invariant under the action.

There are two pleasing consequences of this identification of polygons with points in a Stiefel manifold:

\begin{enumerate}
	\item It is easy to sample (pseudo-)random points from Haar measure on the Stiefel manifold $\St_2(\R^n)$: apply Gram--Schmidt to a pair of $n$-dimensional vectors with entries chosen from a standard Gaussian distribution. This produces a simple algorithm for stochastically generating large ensembles of random polygons. The resulting polygons have edges whose lengths are on the order of $\frac{1}{\sqrt{n}}$ times the diameter of the polygon, which is qualitatively different from polygons produced from, e.g., random Gaussian vertices, which have edge lengths of the same order as the diameter.
	
	\item It is straightforward to find explicit parametrized paths between points in the Stiefel manifold, yielding a procedure for generating efficient interpolations between any two polygons with the same number of edges.
\end{enumerate}

Therefore, the Stiefel manifold is potentially useful both as a generative tool and for finding aesthetically-pleasing morphs between pre-determined shapes.

The \emph{Mathematica} notebook which generated all of the figures in this paper includes reference implementations of all of the algorithms described below and is available at \url{https://shonkwiler.org/bridges2019}.\footnote{There is also a PDF version available at this location for readers who may not have access to \emph{Mathematica}.} It may be helpful to follow along in that notebook to see the implementation details.

\section*{From Polygons to Stiefel Manifolds}

The basic construction is surprisingly straightforward; see~\cite{Cantarella:2017we} for a more detailed but still elementary introduction to these ideas. The goal is a translation-invariant and scale-invariant representation of polygons, so represent an $n$-gon as a sequence of edge vectors $e_1, \ldots , e_n$ which, up to scale, we normalize to have $|e_1| + \ldots + |e_n| = 2$. This representation is really of an \emph{ordered} polygon, since it implies a choice of which edge is first, second, and so on. Notice that a cyclic permutation of the edge labels does not change the visual appearance of the polygon at all.  

In this vector representation, the condition that the polygon closes is equivalent to the vector sum $e_1 + \ldots + e_n$ being the zero vector. Generating random tuples of vectors satisfying $\sum e_k = 0$ and $\sum |e_k| = 2$ is not easy, and finding nice paths in the collection of such tuples is even harder, but a clever trick originally due to Hausmann and Knutson~\cite{Hausmann:1997p8571} (cf.~\cite{Cantarella:2013bla}) comes to the rescue.

Interpret the edge vectors $e_k$ as complex numbers, and then define the complex numbers $z_k$ by $z_k^2 = e_k$. There is ambiguity in the choice of each $z_k$, but the payoff is worth it: letting each $z_k = x_k + \I y_k$ and thinking of $\vec{x} = (x_1, \ldots , x_n)$ and $\vec{y} = (y_1 , \ldots , y_n)$ as vectors in $\R^n$, the closure condition becomes
\[
	0 = \sum e_k = \sum z_k^2 = \sum \left( x_k^2 - y_k^2 + 2\I x_k y_k\right) = \|\vec{x}\|^2 - \|\vec{y}\|^2 + 2 \I \langle \vec{x}, \vec{y}\rangle,
\]
where $\langle \vec{x}, \vec{y}\rangle$ is the standard dot product on $\R^n$, and the normalization condition becomes
\[
	2= \sum |e_k| = \sum |z_k|^2 = \sum \left(x_k^2 + y_k^2\right) = \|\vec{x}\|^2 + \|\vec{y}\|^2.
\]
In other words, a pair $(\vec{x}, \vec{y})$ of $n$-dimensional vectors corresponds to a closed polygon of perimeter 2 if and only if
\[
	\|\vec{x}\| = \|\vec{y}\| = 1 \qquad \text{and} \qquad \langle \vec{x},\vec{y}\rangle = 0.
\]
Pairs of vectors satisfying the above conditions are called \emph{orthonormal 2-frames} for $\R^n$, and the space of all such pairs is the \emph{Stiefel manifold} $\St_2(\R^n)$.

Sampling a random point on the Stiefel manifold $\St_2(\R^n)$ is easy: generate two random vectors $\vec{u}$ and $\vec{v}$ from the standard $n$-dimensional Gaussian distribution and then apply Gram--Schmidt (or another orthogonalization algorithm which preserves spherical symmetry) to get an orthonormal pair $(\vec{x},\vec{y})$. In Mathematica, this can be done with the following function: \begin{verbatim}StiefelSample[n_] := Orthogonalize[RandomVariate[NormalDistribution[], {2, n}]];\end{verbatim}

Alternatively, the same distribution is generated by sampling an $n \times 2$ standard Gaussian matrix $A$ and letting $\vec{x}$ and $\vec{y}$ be the columns of $Q$ in the QR decomposition $A = QR$.\footnote{This requires some care since the QR decomposition is not unique. Due to how NAG and LAPACK implement the QR decomposition~\cite{Lehoucq:1996gh}, Mathematica, Matlab, NumPy, SageMath, Julia, and many other software packages always return a $Q$ with nonpositive $(1,1)$ entry, producing a non-uniform distribution of polygons unless this is corrected for.}

\section*{Random $n$-Gons}

\begin{figure}[h!t]
	\centering
		\includegraphics[height=2in]{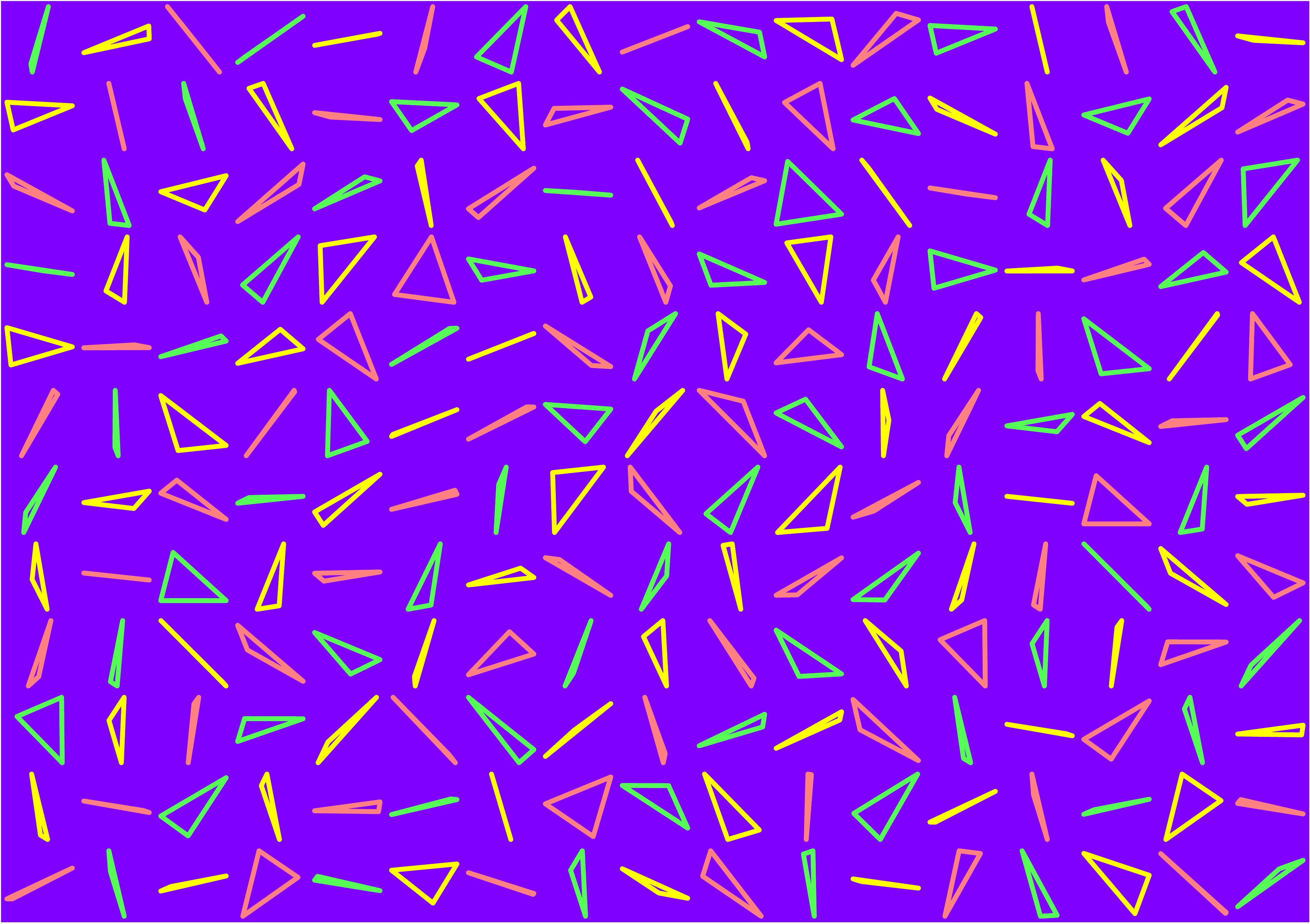}\qquad
	\caption{\emph{Wallpaper} (joint with Anne Harding), showing 204 random triangles, of which 177 (or about 86.8\%) are obtuse, compared to a true probability of $\frac{3}{2}-\frac{3\ln 2}{\pi}\approx 0.838$.}
	\label{fig:triangles}
\end{figure}

Given a point $(\vec{x},\vec{y})$ on the Stiefel manifold, turning it into a polygon is also simple: the $k$th edge of the polygon is simply $(x_k + \I y_k)^2$. Using this technique, it is easy to generate large ensembles of random polygons as in Figure~\ref{fig:triangles}, which illustrates the slightly remarkable fact -- suggested by Wesley Woolhouse~\cite{Woolhouse:1865vw} and Lewis Carroll~\cite[Pillow Problem \#58]{Carroll:1893td} in the 19th century -- that the majority of random triangles are obtuse~\cite{Cantarella:2017we}.

\begin{figure}[h!t]
	\centering
		\includegraphics[height=6in]{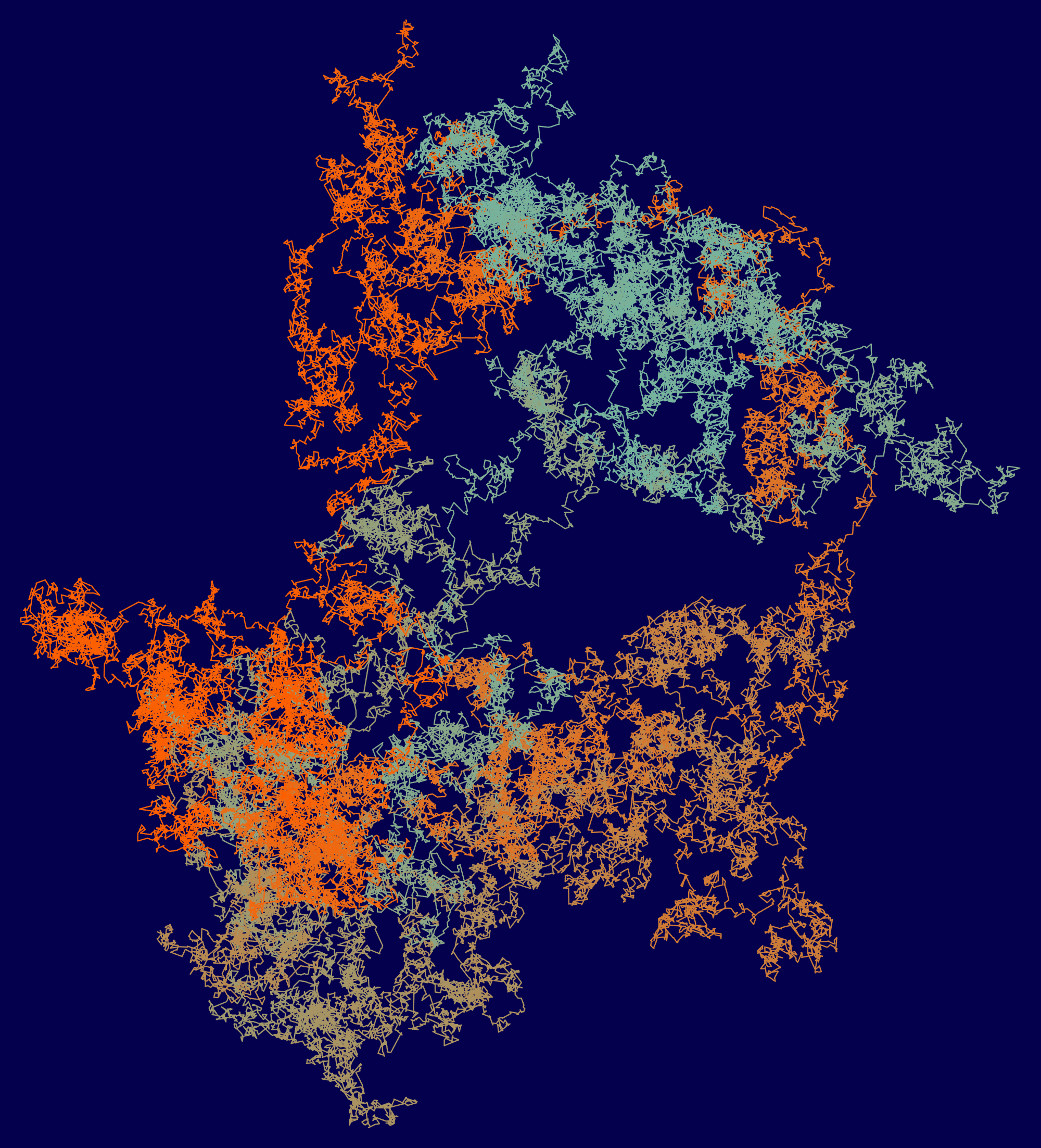}
	\caption{A random 100,000-gon, colored so that the first edge is a bright orange ({\color[RGB]{253,95,0}$\bullet$}), the 50,000th is a light green ({\color[RGB]{118,179,157} $\bullet$}), and the edge colors between are interpolated sinusoidally.}
	\label{fig:100000gon}
\end{figure}

\begin{figure}[h!t]
	\centering
		\includegraphics[height=2.5in]{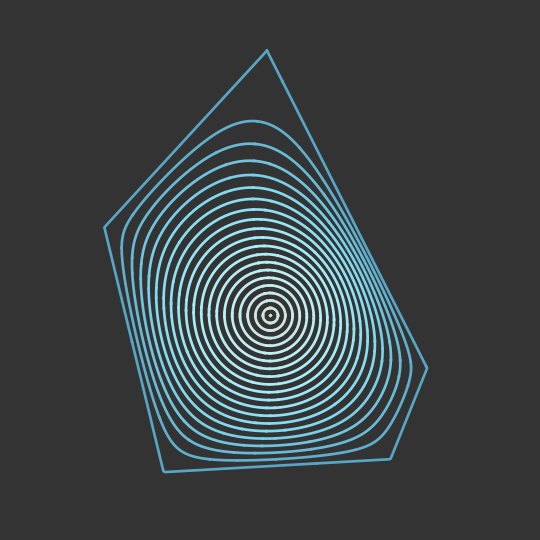}
	\caption{A frame from the animation \emph{Nucleation}~\cite{nucleation}. The Riemann mapping theorem guarantees there is a conformal transformation from the disk to the randomly-generated convex pentagon; the curves show the images of concentric circles under this map.}
	\label{fig:pentagon}
\end{figure}

The 100,000-gon in Figure~\ref{fig:100000gon} illustrates one curious feature of this parametrization of polygon space: the polygons need not be convex or even embedded. However, producing convex $n$-gons is simple: generate a random $n$-gon, then permute the edges until the edge directions are sorted by increasing angle with the $x$-axis. Here is the convexification of a random 200-gon: \includegraphics[height=1em]{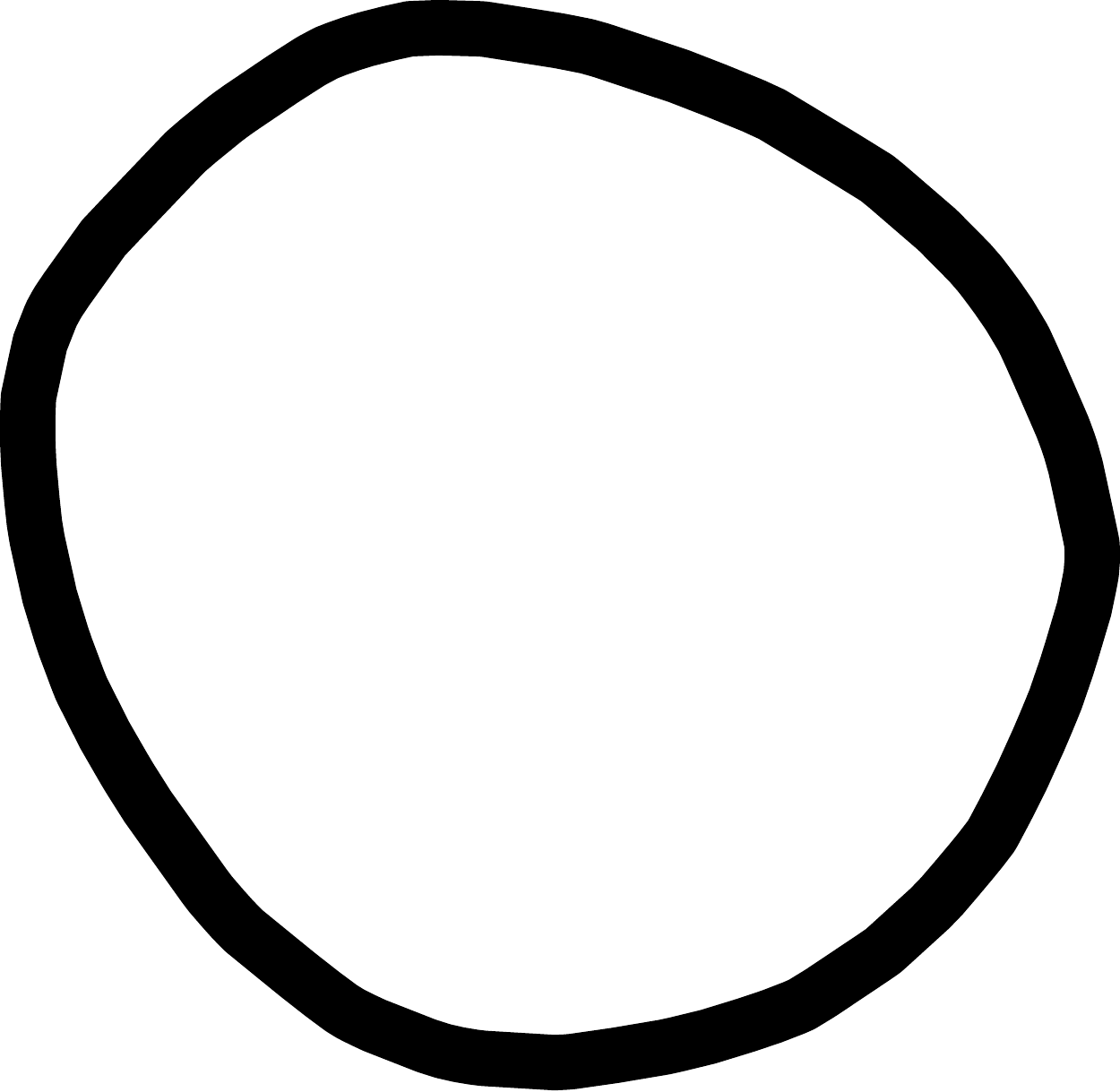}. As discussed in~\cite{Cantarella:2017we}, permuting the edges is an isometry of the Stiefel manifold, so this procedure corresponds to uniformly sampling the convex $n$-gons inside the space of all $n$-gons; see Figure~\ref{fig:pentagon} for an example of a random convex pentagon. Unlike taking the convex hull of a collection of random points, this gives precise control on the number of sides of the resulting convex polygon.

\begin{figure}[h!t]
	\centering
		\includegraphics[height=1.7in]{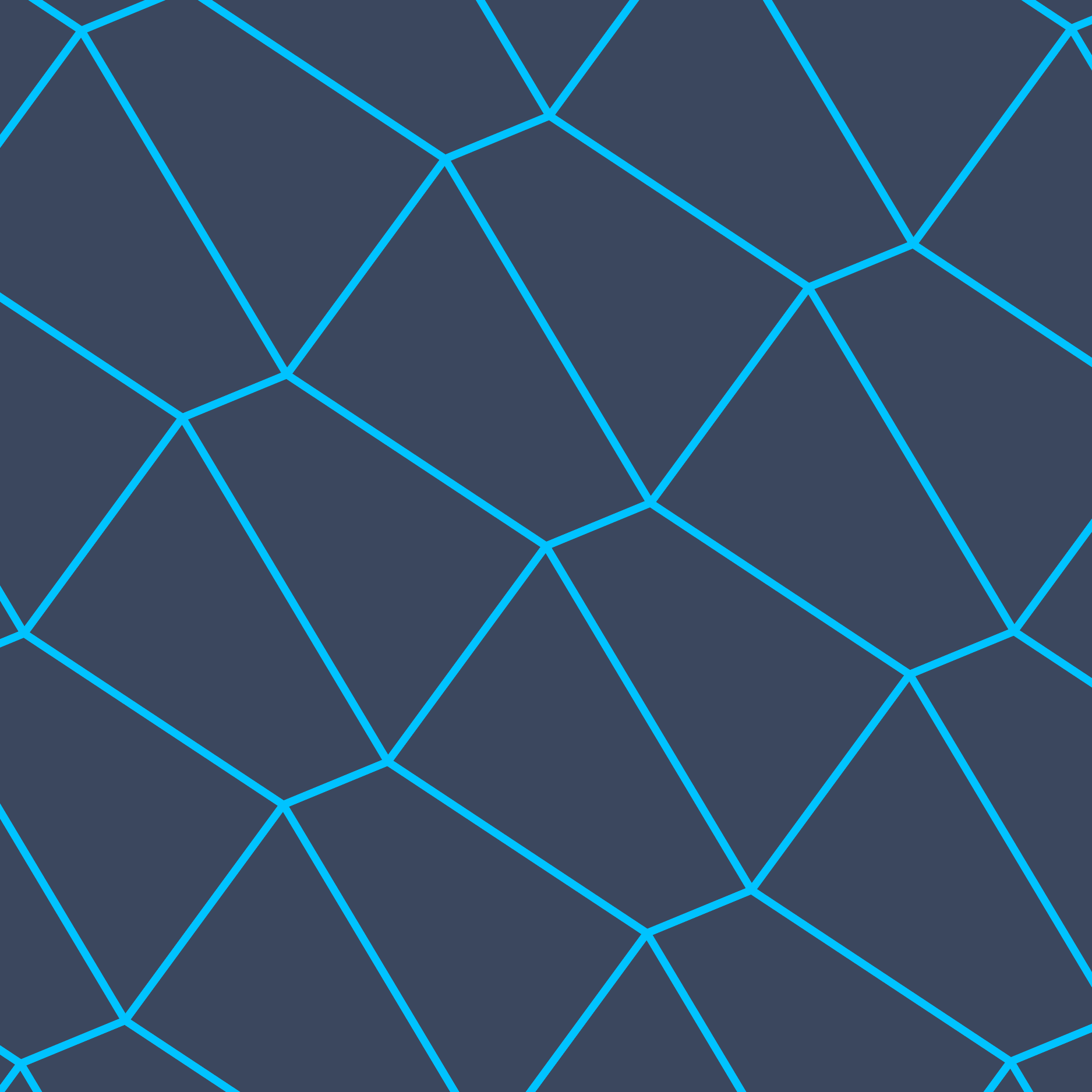}\qquad
		\includegraphics[height=1.7in]{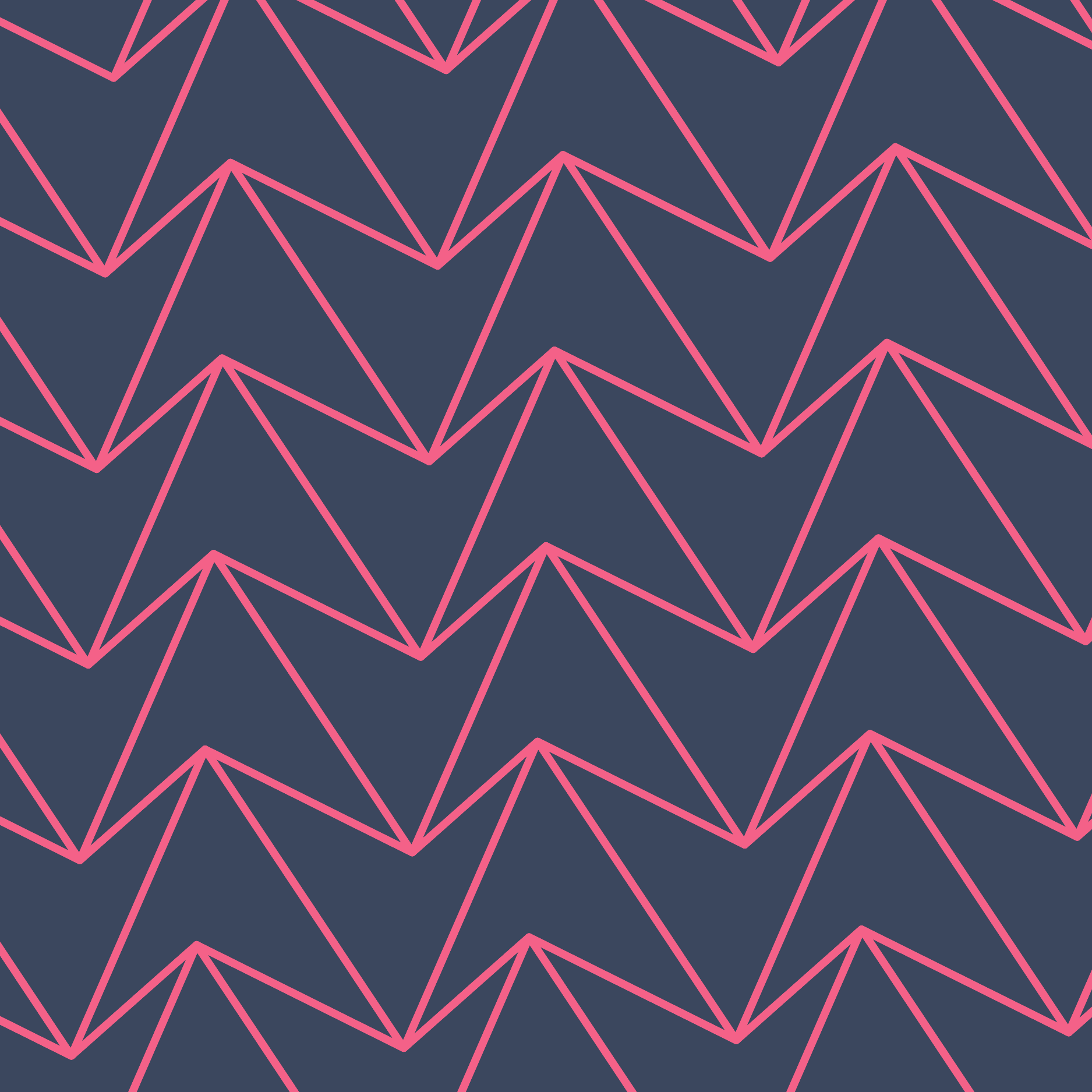}\qquad
		\includegraphics[height=1.7in]{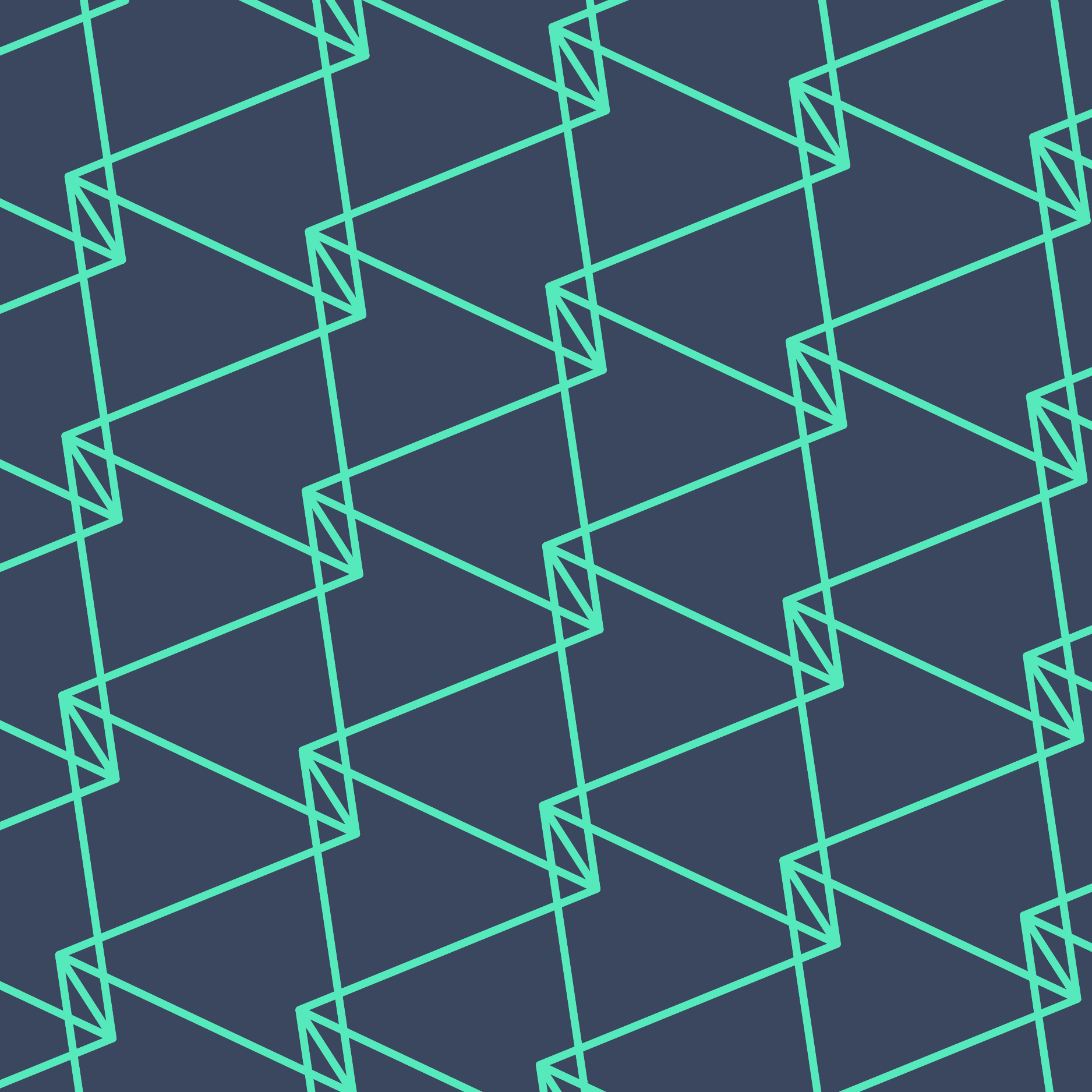}
	\caption{Tilings of the plane by convex, reflex, and crossed quadrilaterals.}
	\label{fig:tilings}
\end{figure}

Since every quadrilateral tiles the plane, the Stiefel manifold approach gives a way to generate random planar tilings: generate a random point in $\St_2(\R^4)$, turn it into a quadrilateral, and then tile the plane with that quadrilateral. Figure~\ref{fig:tilings} shows tilings produced in this way using representatives from each of the three basic classes of quadrilaterals: convex, reflex, and crossed.\footnote{The crossed quadrilaterals do not tile the plane in the traditional sense, as there are usually overlaps. However, they still tile in a generalized sense: it is possible to choose a fixed orientation on the boundary of the quadrilateral and to cover the plane with copies of the quadrilateral-with-oriented-boundary so that at every point in the plane that does not lie on an edge, the sum of the winding numbers of all the quadrilateral boundaries around that point is $+1$.} A beautiful fact proved in~\cite{Cantarella:2017we} (though the result was in some sense anticipated by the science educator, astronomer, future priest, and past Senior Wrangler James Maurice Wilson in 1866~\cite{Wilson:1866vj}) is that, with respect to this measure on quadrilateral space, exactly $\nicefrac{1}{3}$ of quadrilaterals are convex, $\nicefrac{1}{3}$ are reflex, and $\nicefrac{1}{3}$ are crossed.

\section*{Paths in Polygon Space} 
\label{sec:paths_in_polygon_space}

\begin{figure}[h!t]
	\centering
		\includegraphics[height=2in]{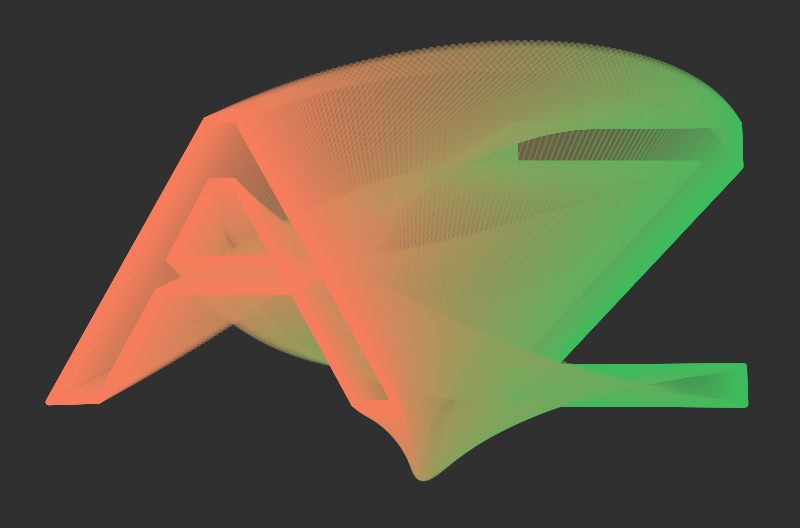} \quad
		\includegraphics[height=2in]{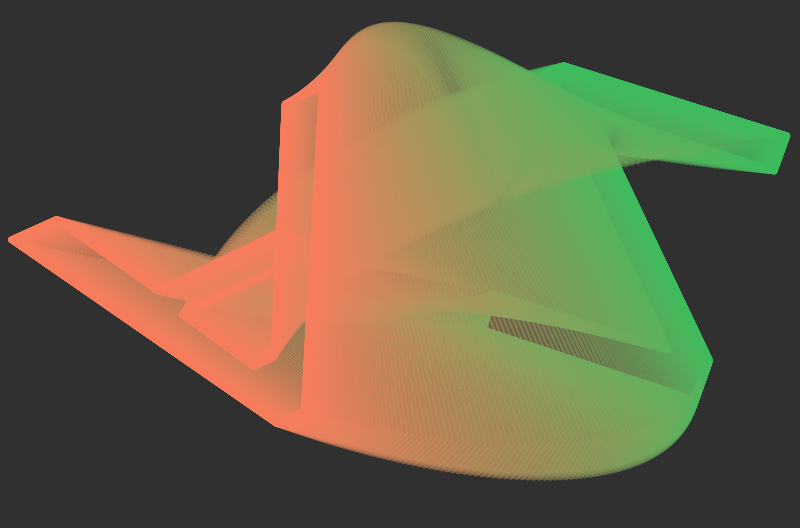}
	\caption{Left: \emph{A to Z}, showing a path in $\St_2(\R^{12})$. Right: the geodesic in $\Gr_2(\R^{12})$ with the same starting and ending polygons, now registered to each other.}
	\label{fig:a to z}
\end{figure}

In addition to ease of sampling, one of the main mathematical motivations for using the Stiefel manifold model for polygons is that it is straightforward to find paths between pairs of points on the Stiefel manifold, and hence between any two given $n$-gons. An obvious choice for such a path would be a geodesic, but unfortunately it is not generally possible to explicitly parametrize geodesics between arbitrary pairs of points in Stiefel manifolds. Instead, the following procedure, used to produce Figure~\ref{fig:a to z}, is modeled off the construction of geodesics in the Grassmannian discussed below, and hence is a reasonably efficient strategy for traveling between two points, though it almost never produces a geodesic. Given the starting and ending points $(\vec{x}_0,\vec{y}_0),(\vec{x}_1, \vec{y}_1) \in \St_2(\R^n)$, let $\vec{x}(t)$ be the direct rotation of $\vec{x}_0$ towards $\vec{x}_1$ and let $\vec{y}(t)$ be the normalized projection of the direct rotation of $\vec{y}_0$ towards $\vec{y}_1$ onto the perpendicular space $\vec{x}(t)^\bot$. Then $(\vec{x}(t),\vec{y}(t))$ is a path in $\St_2(\R^n)$ from $(\vec{x}_0,\vec{y}_0)$ to $(\vec{x}_1, \vec{y}_1)$.\footnote{Reversing the roles of the $\vec{x}_i$ and the $\vec{y}_i$ in this procedure generally produces a very similar path.}

A more mathematically satisfying approach is to pass from $\St_2(\R^n)$ to the corresponding \emph{Grassmann manifold} $\Gr_2(\R^n)$ of 2-dimensional linear subspaces of $\R^n$, sending a pair of orthonormal vectors to the subspace they span. This gives an orientation-independent representation of polygons~\cite{Cantarella:2013bla}, meaning that it does not distinguish between a polygon and a rotated version of the same polygon. 

Geodesics in the Grassmannian are well-understood~\cite{Bjorck:1973jq,Edelman:1999ei}: if ${P_0, P_1 \in \Gr_2(\R^n)}$, $M_i$ is an $n \times 2$ matrix with orthonormal columns spanning $P_i$, and $U \Sigma V^T$ is the singular value decomposition of $M_0^T M_1$, let $(\vec{x}_0,\vec{y}_0)$ be the columns of $M_0 U$ and $(\vec{x}_1,\vec{y}_1)$ be the columns of $M_1 V$. Then $(\vec{x_i},\vec{y}_i)$ is an orthonormal basis for $P_i$, and the geodesic in $\Gr_2(\R^n)$ from $P_0$ to $P_1$ corresponds to the path $(\vec{x}(t),\vec{y}(t))$ in $\St_2(\R^n)$ where $\vec{x}(t)$ is the direct rotation of $\vec{x}_0$ towards $\vec{x}_1$ and $\vec{y}(t)$ is the direct rotation of $\vec{y}_0$ towards $\vec{y}_1$. The effect of choosing a different orthonormal basis for an element of $G_2(\R^n)$ is to rotate the corresponding polygon; the bases for the starting and ending points that arise as above optimally register the polygons to each other. This loss of control of orientation may not be desirable if the original orientations were important.

One curious feature of the Stiefel manifold parametrization of polygon space -- which is either a detriment or an opportunity, depending on circumstance -- is that there are many choices in how to lift a polygon to a point in the Stiefel manifold. Since this is really a parametrization of \emph{ordered} polygons, it depends on an ordering of the edges. This means that cyclically permuting the edge labels produces $n$ different paths between visually identical starting and ending polygons.

\begin{figure}[h!t]
	\centering
		\includegraphics[height=2in]{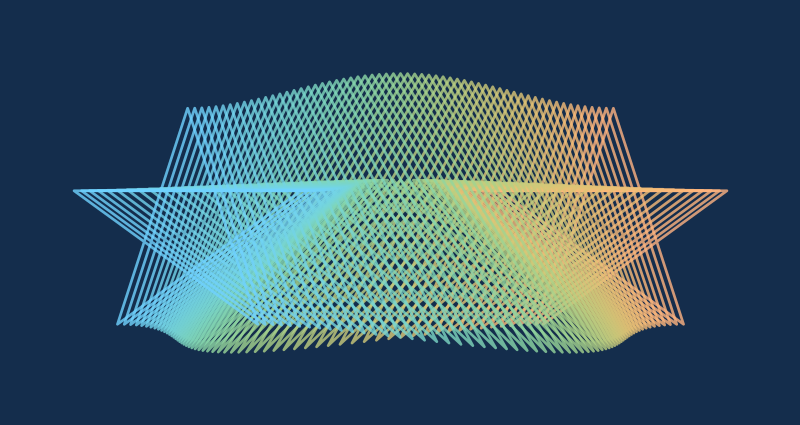}
	\caption{The geodesic in $\Gr_2(\R^5)$ between two different lifts of the same pentagonal star.}
	\label{fig:star to star}
\end{figure}

In fact, there are even more choices: lifting a polygon involved interpreting each edge vector $e_k$ as a complex number and finding $z_k$ so that $z_k^2 = e_k$. There are generally two possible choices for the square root, so this produces $2^n$ possible lifts of a given ordered $n$-gon. One consequence is that there are many non-trivial paths from a polygon to itself coming from paths between different lifts. Figure~\ref{fig:star to star} shows such a path between different lifts of a pentagonal star.

\section*{A Generalization to $n$-Gons in Space} 
\label{sec:a_generalization_to_n_gons_in_space}

A generalization of this construction yields a correspondence between polygons in $\R^3$ and points in the Stiefel manifold~$\St_2(\C^n)$. The idea is essentially the same as for planar polygons, except now instead of thinking of the edge vectors as complex numbers, we think of them as purely imaginary quaternions, and instead of taking a square root we take a section of the Hopf map $q \mapsto \overline{q}\,\I\, q$. In this way a closed polygon in space can be turned into a pair $(\vec{x},\vec{y})$ of vectors in $\C^n$ which are orthonormal with respect to the standard Hermitian inner product on $\C^n$. Specifically, $(\vec{x},\vec{y}) \in \St_2(\C^n)$ is mapped to the polygon in $\R^3$ whose $\ell$th edge is $(\overline{x}_\ell - y_\ell \J) \I (x_\ell+y_\ell \J)$, which is a purely imaginary quaternion, meaning it can be written as $a \I + b \J + c \K$ for some $a$, $b$, and $c$, and hence interpreted as a three-dimensional vector.\footnote{Here $\I, \J$, and $\K = \I\J$ are the standard quaternion units.}  See~\cite{Cantarella:2013bla} for details.

\begin{figure}[h!t]
	\centering
		\includegraphics[height=1.275in]{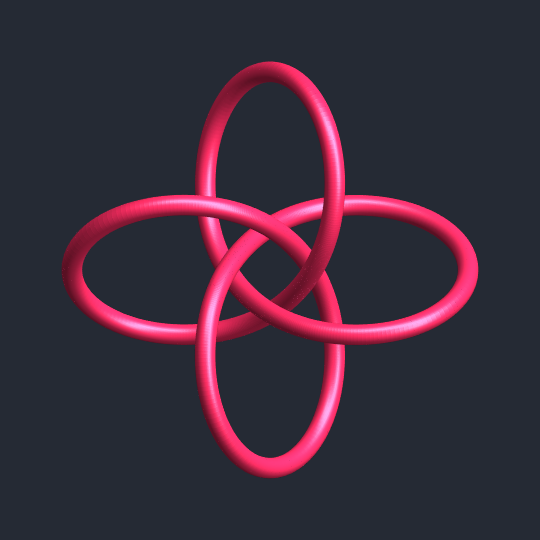}
		\includegraphics[height=1.275in]{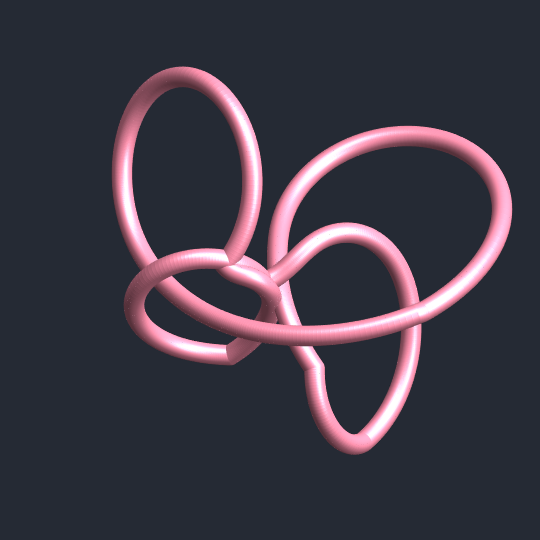}
		\includegraphics[height=1.275in]{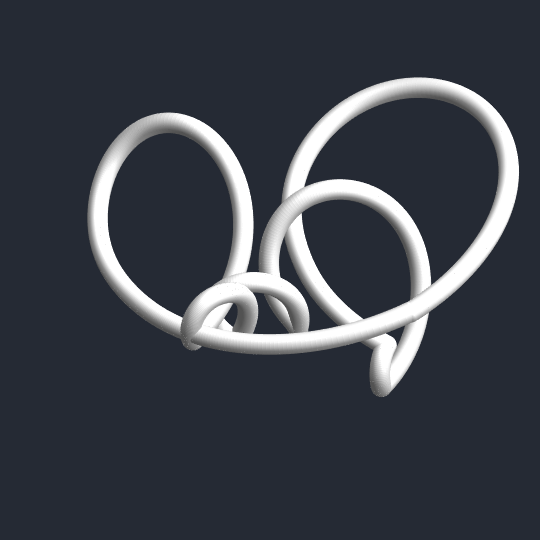}
		\includegraphics[height=1.275in]{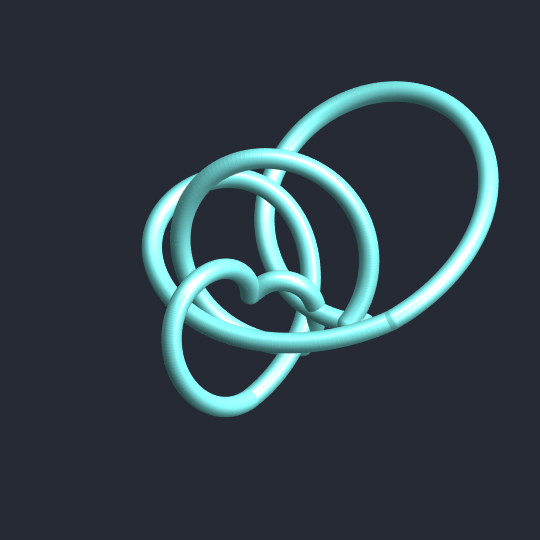}
		\includegraphics[height=1.275in]{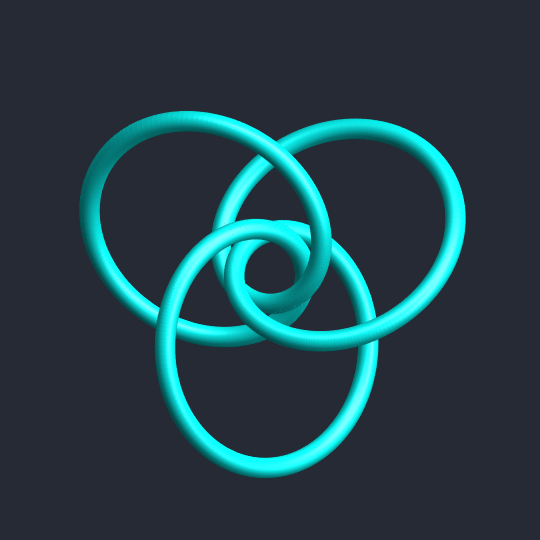}
	\caption{A path in $\St_2(\C^{1000})$ from a $(3,4)$-torus knot to a $(4,3)$-torus knot.}
	\label{fig:morph}
\end{figure}

Again, we can both sample random $n$-gons and, since it is straightforward to find paths connecting pairs of points in $\St_2(\C^n)$ or geodesics in the corresponding Grassmannian $\Gr_2(\C^n)$, morph any polygon in $\R^3$ to any other in a natural way, as in Figure~\ref{fig:morph}, which shows a path in $\St_2(\C^{1000})$ between 1000-gons which form $(3,4)$ and $(4,3)$ torus knots.

\begin{figure}[h!t]
	\centering
		\includegraphics[height=4in]{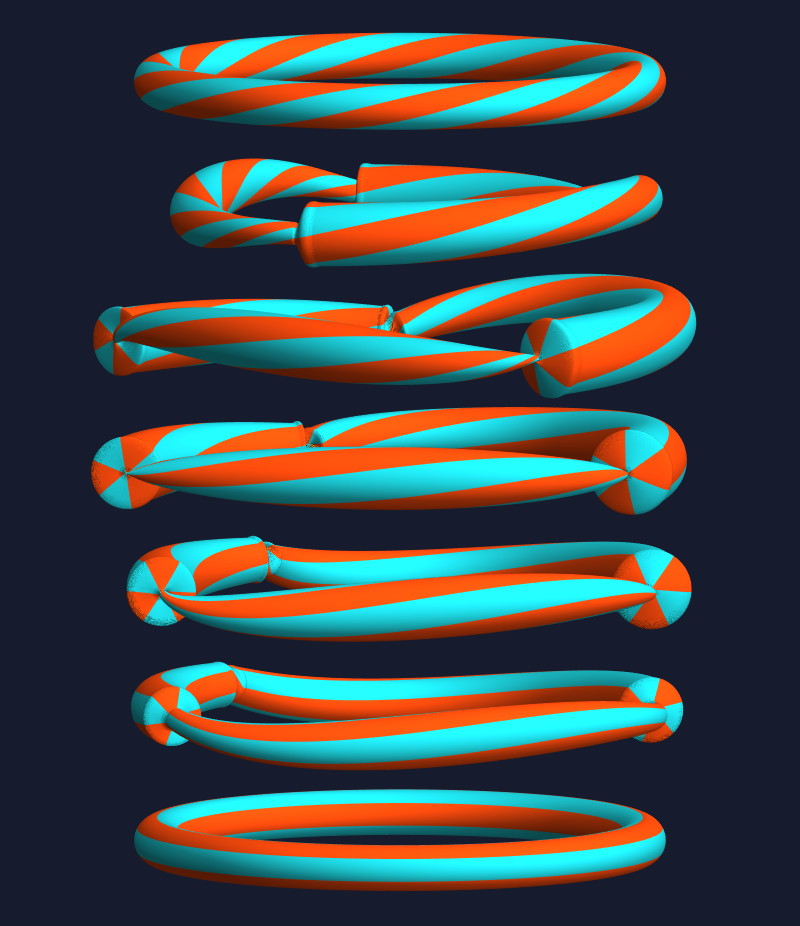}
	\caption{The geodesic in the Grassmannian $\Gr_2(\C^{200})$ between the regular 200-gon with the standard blackboard framing (bottom) and the regular 200-gon with a thrice-twisted frame (top).}
	\label{fig:framing geodesic}
\end{figure}

Just as the lift of a planar polygon depended on choices of square roots, lifting a space polygon to~$\St_2(\C^n)$ depends on choices of sections of the Hopf map, meaning there is an entire circle's worth of possible preimages for each edge. This circle corresponds to possible framings on each edge, and lifts of polygons can be interpreted as \emph{inflatable elastic knots}, which are bendable, twistable closed curves with variable thickness. Geodesics in the Grassmannian minimize an appropriate generalized Kirchhoff energy~\cite{Needham:2017vd} which penalizes bending, twisting, and varying the thickness. Figure~\ref{fig:framing geodesic} shows a geodesic between two different framings of the regular 200-gon: the apparent singularities arise because it is energetically advantageous to concentrate most of the bending and twisting in a few locations and moreover because the energy being minimized penalizes bending and twisting more where the knot is thicker, so these locations tend to have (nearly) zero cross-sectional radius. Notice, in particular, that there are three singularities, corresponding to the fact that the ending frame has three twists in it. 

Natural choices of framing for a polygonal curve include the discrete Frenet--Serret frame (which often approximates a material frame, for example in proteins~\cite{Hu:2011bxa}) and the discrete Bishop frame (which minimizes twist), but the framing ambiguity provides an opportunity for experimentation.


\section*{Summary and Conclusions}

The Stiefel manifold parametrizations of polygons both in the plane and in space provide new tools for generative production of polygons and for morphing one polygon into another. While the primary motivation for developing this machinery was to find an exactly solvable model for ring polymers (see~\cite{Cantarella:2013bla} for much more detail and background), I believe that this approach should also be useful for both generative and deterministic art, especially since it requires at most a bit of linear algebra to implement.

\section*{Acknowledgements}

I am very grateful to the three anonymous referees whose thoughtful comments have substantially improved this paper. This work was supported by a grant from the Simons Foundation (\#354225).

{\setlength{\baselineskip}{13pt} 
\raggedright				

} 

\end{document}